\documentclass{article}

\title{The farthest point map on the $4$-cube}
\author{Yoshikazu Yamagishi\footnote{yg@rins.ryukoku.ac.jp, 
Faculty of Advanced
Science and Technology, Ryukoku University,
Seta, Otsu, 520-2194, Japan}}

\usepackage{amsfonts}
\usepackage{amssymb}
\usepackage{amsmath}
\usepackage{amsthm}
\usepackage[dvipdf]{graphicx}

\newtheorem{theorem}{Theorem}

\newtheorem{lemma}[theorem]{Lemma}
\newtheorem{proposition}[theorem]{Proposition}

\newcommand{\BackRight}{{\rm BackRight}}
\newcommand{\Back}{{\rm Back}}
\newcommand{\Front}{{\rm Front}}

\newcommand{\aff}{{\rm Aff}}

\newcommand{\dist}{{\rm d}}

\newcommand{\half}{\frac{1}{2}}

\newcommand{\id}{{\rm id}}

\newcommand{\incircle}{{\rm incircle}}

\newcommand{\insphere}{{\rm insphere}}

\newcommand{\lbl}{{\rm label}}

\newcommand{\outcircle}{{\rm outcircle}}
\newcommand{\outsphere}{{\rm outsphere}}
\newcommand{\pBD}{p_{\sfB\sfD}}

\newcommand{\pB}{p_\sfB}
\newcommand{\pD}{p_\sfD}
\newcommand{\pFD}{p_{\sfF\sfD}}
\newcommand{\pF}{p_\sfF}
\newcommand{\pLD}{p_{\sfL\sfD}}

\newcommand{\pL}{p_\sfL}
\newcommand{\pRD}{p_{\sfR\sfD}}

\newcommand{\pR}{p_\sfR}

\newcommand{\pUB}{p_{\sfU\sfB}}
\newcommand{\pUF}{p_{\sfU\sfF}}

\newcommand{\pUL}{p_{\sfU\sfL}}

\newcommand{\pUR}{p_{\sfU\sfR}}
\newcommand{\pU}{p_\sfU}
\newcommand{\psird}{{\psi_1}}
\newcommand{\psiur}{{\psi_2}}

\newcommand{\qBD}{q_{\sfB\sfD}}

\newcommand{\qB}{q_\sfB}
\newcommand{\qD}{q_\sfD}

\newcommand{\qLD}{q_{\sfL\sfD}}

\newcommand{\qL}{q_\sfL}
\newcommand{\qRD}{q_{\sfR\sfD}}

\newcommand{\qR}{q_\sfR}

\newcommand{\qUB}{q_{\sfU\sfB}}

\newcommand{\qUL}{q_{\sfU\sfL}}

\newcommand{\qUR}{q_{\sfU\sfR}}
\newcommand{\qU}{q_\sfU}

\newcommand{\real}{{\mathbb R}}

\newcommand{\rmc}{{\rm c}}

\newcommand{\sBD}{{\sf BD}}
\newcommand{\sBLD}{{\sf BLD}}
\newcommand{\sBL}{{\sf BL}}
\newcommand{\sBRD}{{\sf BRD}}
\newcommand{\sBR}{{\sf BR}}
\newcommand{\sB}{{\sf B}}
\newcommand{\sD}{{\sf D}}
\newcommand{\sFD}{{\sf FD}}
\newcommand{\sF}{{\sf F}}
\newcommand{\sG}{{\sf G}}
\newcommand{\sLD}{{\sf LD}}
\newcommand{\sLFD}{{\sf LFD}}
\newcommand{\sLF}{{\sf LF}}
\newcommand{\sL}{{\sf L}}
\newcommand{\sRD}{{\sf RD}}
\newcommand{\sRFD}{{\sf RFD}}
\newcommand{\sRF}{{\sf RF}}
\newcommand{\sR}{{\sf R}}
\newcommand{\sS}{{\sf S}}
\newcommand{\sUBL}{{\sf UBL}}
\newcommand{\sUBR}{{\sf UBR}}
\newcommand{\sUB}{{\sf UB}}
\newcommand{\sUF}{{\sf UF}}
\newcommand{\sULF}{{\sf ULF}}
\newcommand{\sUL}{{\sf UL}}
\newcommand{\sURF}{{\sf URF}}
\newcommand{\sUR}{{\sf UR}}
\newcommand{\sU}{{\sf U}}
\newcommand{\set}[1]{\{ #1 \}}
\newcommand{\sfBD}{{\sf BD}}
\newcommand{\sfB}{{\sf B}}
\newcommand{\sfD}{{\sf D}}
\newcommand{\sfF}{{\sf F}}
\newcommand{\sfG}{{\sf G}}

\newcommand{\sfL}{{\sf L}}
\newcommand{\sfRD}{{\sf RD}}
\newcommand{\sfR}{{\sf R}}
\newcommand{\sfS}{{\sf S}}

\newcommand{\sfUB}{{\sf UB}}

\newcommand{\sfUR}{{\sf UR}}
\newcommand{\sfU}{{\sf U}}

\newcommand{\src}{{\rm src}}

\newcommand{\radius}{{\rm rad}}
\newcommand{\diam}{{\rm diam}}

\begin{document}

\maketitle

\begin{abstract}
We study the farthest point mapping on (the boundary of) the $4$-cube
with respect to the intrinsic metric, and its dynamics as a
multivalued mapping.  It is a piecewise rational map.  It is more
complicated than the one on the $3$-cube, but it is shown that the limit
set of the farthest point map on the $4$-cube is the union of the
diagonals of eight ($3$-cube) facets, like the farthest point map on the
$3$-cube whose limit set is the union of the six (square) facets. This
is in contrast to the doubly covered simplices and (the boundary of)
the regular $4$-simplex, where the limit set is a finite set.  If the
source point is in the interior of a facet, its limit set is also in
the facet.

The farthest point mapping is closely related to the star unfolding 
and source unfolding. We give a loose definition
of star unfolding of the surface of a $4$-dimensional polytope.

We also study the intrinsic radius and diameter of the $4$-cube.  It
is expected that the intrinsic radius/diameter ratio of an $n$-cube is
monotonically decreasing in dimension.
\end{abstract}

\section{Introduction}

This paper studies the farthest point map on the $4$-cube $I^4$, where $I=[0,1]$.
It is natural to ask for the farthest point in a compact metric space.
A convex surface is called a Steinhaus surface if the farthest point map is single-valued and involutive \cite{croft,rouyer-steinhaus}.  
Examples of a Steinhaus surface other than the sphere
was given by \cite{itoh-antipodal,itoh-cylinders,vilcu-steinhaus,vilcu-symmetry}.
Rouyer showed that no convex polytope is Steinhaus \cite{rouyer-steinhaus}.
Thus, the dynamics of the farthest point map on convex polytopes are of interest.

There are not many concrete examples of the farthest point map 
on convex polytopes, especially in higher dimensions.
Farthest point map on the regular tetrahedron is 
studied in \cite{rouyer-tetrahedre},
flat surfaces in \cite{rouyer-flat},
regular octahedron in \cite{schwartz-octa},
regular dodecahedron in \cite{schwartz-dodeca},
centrally symmetric polyhedron in \cite{wang},
doubly covered simplex and convex polyhedron in \cite{itoh-degenerate},
and doubly covered parallelotope in \cite{ueda-yg}.

Let $M = \partial P$ be the surface of a convex polytope $P \subset \real^{n+1}$.
The intrinsic metric $\dist$ on $M$ is defined by the length of 
the shortest paths.
Let $p\in M$, and call it a source point.
The set of the farthest points of $p$ is defined by
\[ f(p) = \set{q \in M \mid 
\dist(p,q) \ge \dist(p,q'), \forall q' \in M}.
\]
An orbit of $p$ is defined as a sequence $\set{q_i}_{i\ge0}$
by $q_0=p$, $q_{i+1}\in f(q_i)$, $i\ge0$.

The farthest point map is closely related to the cut locus.
A point $q \in M$ is called a cut point if any shortest path from $p$ to $q$ 
cannot be extended as a shortest path from $p$.
A farthest point is a cut point.
The set of cut points is called the cut locus of $p$ and is denoted by $C(p)$.
It is the closure of the set of points that have more
than one shortest path from $p$.
A point in an $(n-2)$-face of $M$ is called a warped point. 
A warped point is a cut point \cite{miller-pak}.
A recent result on the length of the cut locus on convex surfaces is
given in \cite{yuan-zamfirescu}.

The farthest point map is also related to star unfolding and source unfolding \cite{demaine-orourke-book,miller-pak}.  
The computation of the farthest point map consists of several steps as follows.
First, the source images are given in the affine hull of each facet
as vertices of the star unfolding.
Second, the Voronoi tessellation with respect to the source images gives the cut locus.
Third, the cut locus gives the source unfolding and the farthest point map.
The star unfolding has not been defined in higher dimensions $n\ge3$.
In this paper we try to define the star unfolding in higher dimensions.

The source unfolding of $M$ is defined as follows.
Let $M_1\subset \real^n$ be a (not convex) polytope. 
Let $g_1 : M_1 \to M$ be a surjective continuous map
 which is piecewise isometry.
Denote the interior of $M_1$ by $M_1^o$. 
Suppose that 
$g_1|M_1^o$ is a homeomorphism of $M_1^o$ onto $M \setminus C(p)$.
Then $M_1$ (or more precisely, the inverse mapping
$(g_1|M_1^o)^{-1} : M \setminus C(p) \to M_1^o$) is called a source unfolding of $M$ with respect to $p$.
The sphere in $\real^n$ with center $g_1^{-1}(p)$ 
and radius $\dist(p,f(p))$ circumscribes $M_1$,
and a contact point $q' \in \partial M_1$ corresponds to a farthest point $g_1(q') \in f(p)$.

The notion of star unfolding of a convex polyhedron is defined when $n=2$,
but not when $n\ge3$.  
In this article we try to give a loose definition of star unfolding as follows. 
Let $M_2 \subset \real^n$ be a polytope.
Let $g_2 : M_2 \to M$ be a surjective continuous map
which is piecewise isometry.
Suppose that 
$g_2|M_2^o : M_2^o \to M \setminus S$ is a homeomorphism,
where $S$ is a union of some (infinitely many) shortest paths with an endpoint $p$.
Then $M_2$ (or more precisely, the inverse mapping 
$(g_2|M_2^o)^{-1} : M \setminus S \to M_2^o$) is called a star unfolding of $M$ with respect to $p$.
The points in $g_2^{-1}(p)$ are called the source images.
The star unfolding $M_2$ is called maximal
 if for any star unfolding
$g'_2: M'_2 \to M$ with a homeomorphism 
$g'_2|(M'_2)^o : (M'_2)^o \to M\setminus S'$ such that $S' \subset S$, we have
$S'=S$.

If $n=2$, $S$ is defined as the union of the shortest paths from $p$ to the vertices of $M$.  
If $n=3$, the set of warped points is
the the union of the $1$-faces of $M$.
It is not enough for $S$ to be the union of the shortest paths from $p$ to the warped points,
because there are some warped points $q$ that have more than one shortest path from $p$.
We need additional cuts between these shortest paths
as shown in Figure~\ref{fig:4cube-star},
which may not be unique.

This paper describes the farthest point map on the $4$-cube.
It is 
shown that the limit set of the farthest point map is the union of the 
diagonals of the eight ($3$-cube) facets.
Section 2 recalls the star and source unfolding and 
the farthest point map on the $3$-cube.  
The limit set of the farthest point map on the $3$-cube is the union
of the diagonals.
Section 3 gives the $26=3^3-1$ source images in the affine hull of
the `goal' facet, and a star and source unfolding of the $4$-cube.
Sections 4 and 5 show that it is sufficient to consider 
only $11$ out of the $26$ source images 
as the site points whose Voronoi domains are adjacent to
the farthest points.
In Section 6 we show that all the farthest points are adjacent to
the Voronoi domain of the site called $\pU$.
In Section 7 we determine the farthest 
points by comparing the distance from $\pU$.
The farthest point map is a piecewise continuous map 
consisting of $10$ rational maps.
In Section 8 we show that 
the domain for the rational map denoted by $\rmc(\sB,\sD,\sR)$
is forward invariant under the dynamics.
Its limit set is the diagonal
of the ($3$-cube) facet.
For any source point $p$ inside the facet,
any orbit falls into the domain of $\rmc(\sB,\sD,\sR)$,
which implies that the limit set of the farthest point map is the
union of the diagonals.
It is also shown that if $p$ is inside the facet, the limit set is
also inside the facet.

Finally, in Section 9 we study the intrinsic radius and 
diameter of the $4$-cube.
It is expected that the intrinsic radius/diameter ratio of the $n$-cube is monotonically decreasing in dimension.

One of the important problems in origami mathematics is to unfold polyhedra without overlapping.
The star unfolding has no overlap when $n=2$.
Miller and Pak \cite{miller-pak} showed that
the source unfolding has no overlap.
It is not proved that the star unfolding has no overlap when $n\ge3$.

Here we give some technical preparation.
The dual of the Voronoi tessellation is called Delaunay tetrahedralization (or triangulation) \cite{cheng-dey-shewchuk}.
In general, let $\Gamma \subset \real^n$ be a site set.
A simplex with vertices $v_1,\dots,v_k \in \Gamma$, $k \le n$,
 is called (weakly) Delaunay
if it has the empty circle property. 
That is, there exists a sphere $C$
such that $v_1,\dots, v_k \in \partial C$
and all other sites $v \in \Gamma$ are outside or on the boundary of $C$.

Denote the affine hull of a subset $L \subset \real^{n+1}$ by $\aff(L)$.
For two adjacent facets $F,F'$ of $M$, 
an (un)folding map 
$\phi_{FF'} : \aff(F') \to \aff(F)$
 is defined by the following conditions \cite{miller-pak},
\begin{itemize}
\item $\phi_{FF'}$ is an isometry,
\item $\phi_{FF'}|(F \cap F') = \id$,
\item $\phi_{FF'}(F') \neq F$.
\end{itemize}
For a sequence of facets $L = (F_1, F_2, \dots, F_\ell)$,
the unfolding map of $\aff(F_\ell)$ onto $\aff(F_1)$
is defined by
\[ \phi_L = \phi_{F1,F2} \circ \dots\circ \phi_{F_{\ell-1},F_\ell} 
: \aff(F_\ell) \to \aff(F_1) .
\]

For $p_i=(x_i,y_i,1)$, $i=1,\dots,4$, let 
\begin{align*}
\outcircle(p_1, p_2,p_3, p_4)
&= 
\begin{vmatrix}
 x_1 & y_1 & 1 & x_1^2 + y_1^2 \\
 x_2 & y_2 & 1 & x_2^2 + y_2^2 \\
 x_3 & y_3 & 1 & x_3^2 + y_3^2 \\
 x_4 & y_4 & 1 & x_4^2 + y_4^2 
\end{vmatrix},
\\
\incircle(p_1,p_2,p_3,p_4) &= - \outcircle(p_1,p_2,p_3,p_4).
\end{align*}
If $\det(p_1,p_2,p_3)
\cdot 
\incircle(p_1,p_2,p_3,p_4)
>0$, 
then
$p_4$ is inside the circumcircle of $p_1,p_2,p_3$,
\cite{cheng-dey-shewchuk}.

For $p_i=(x_i,y_i,z_i,1)$, $i=1,\dots, 5$, let
\begin{align*}
\outsphere(p_1,p_2,p_3,p_4,p_5)
&= 
\begin{vmatrix}
 x_1 & y_1 & z_1 & 1 & x_1^2+y_1^2+z_1^2 \\
 x_2 & y_2 & z_2 & 1 & x_2^2+y_2^2+z_2^2 \\
 x_3 & y_3 & z_3 & 1 & x_3^2+y_3^2+z_3^2 \\
 x_4 & y_4 & z_4 & 1 & x_4^2+y_4^2+z_4^2 \\
 x_5 & y_5 & z_5 & 1 & x_5^2+y_5^2+z_5^2  
\end{vmatrix},
\\
\insphere(p_1,p_2,p_3,p_4,p_5) &= - \outsphere(p_1,p_2,p_3,p_4,p_5).
\end{align*}
If $\det(p_1,p_2,p_3,p_4)
\cdot \insphere(p_1,p_2,p_3,p_4,p_5)>0$, then
$p_5$ is inside the circumsphere of $p_1,p_2,p_3,p_4$.

For a hyperplane $H \subset \real^{n+1}$, 
let $r_H: \real^{n+1} \to \real^{n+1}$
 be the reflection with respect to $H$.

\begin{lemma}
\label{3}
Let $M = \partial P$ be the surface of a convex polytope $P \subset \real^{n+1}$.
Suppose that $M$ is symmetric about a hyperplane $H$.
Let $p,q\in M$, and suppose that 
$q$ is a farthest poit of $p$.
If $p,q\in M\setminus H$, then $q$ is on the opposite side of $p$ with respect to $H$. 
If $p \in H$, then $r_H(q)$ is also a farthest point of $p$.
\end{lemma}
\begin{proof}
Suppose that $q$ is on the same side of $p$ with respect to $H$.
Let $\gamma : [0,1] \to M$ be a shortest path from 
$\gamma(0) = p$ to $\gamma(1) = r_H(q)$.
There exists $t_0 \in (0,1)$ such that $\gamma(t_0) \in H$,
and such that $\gamma(t_0-\epsilon), \gamma(t_0+\epsilon)$ are on the opposite sides with respect to $H$ for small $\epsilon>0$.
Define a path $\tilde\gamma$ by $\tilde\gamma(t) = \gamma(t)$ for $0\le t \le t_0$, and $\tilde\gamma(t) = r_H(\gamma(t))$ for $t_0 \le t \le 1$.
Then $\tilde\gamma$ is a shortest path from $p$ to $q$.
However, 
$\tilde\gamma(t-\epsilon)$, 
$\tilde\gamma(t+\epsilon)$ are on the same side with respect to $H$, so we get a shorter path by connecting
$\tilde\gamma(t-\epsilon)$ and $\tilde\gamma(t+\epsilon)$ by a line segment, a contradiction.

The proof for $p\in H$ is obvious.
\end{proof}

\section{Star and source unfolding and farthest point map on the $3$-cube}

\begin{figure}[t]
  \centering
    \includegraphics[width=40mm]{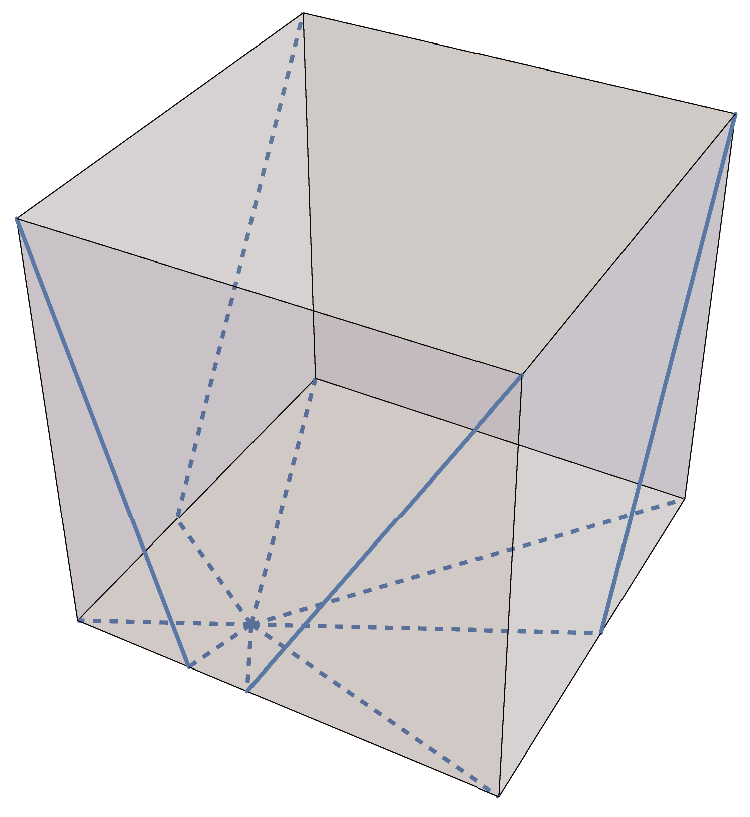}
  \includegraphics[width=40mm]{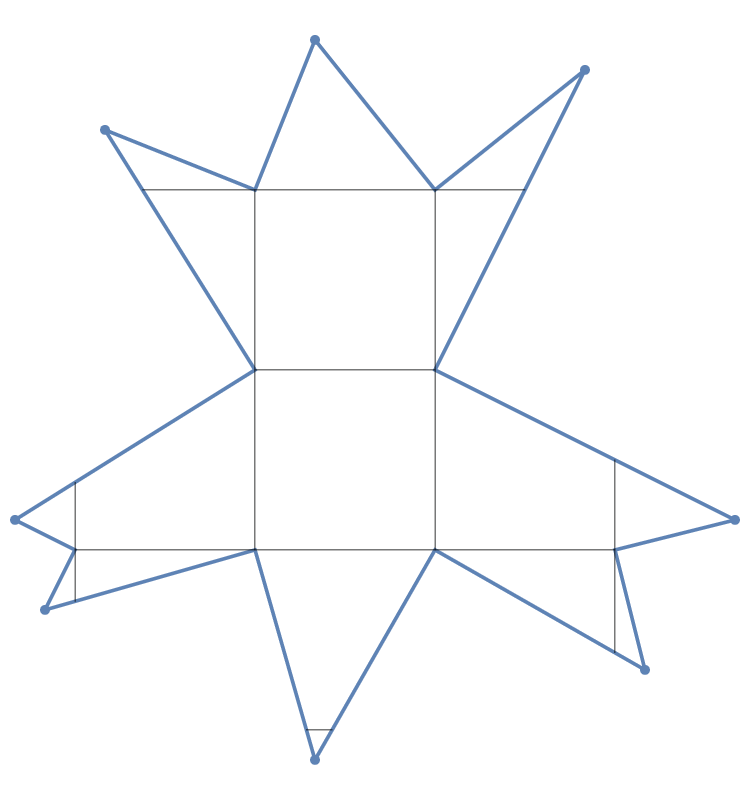}
  \caption{(a) Shortest paths from a source point $p=(1/3,1/6,0)$ to the vertices of the $3$-cube.
  (b) A star unfolding of the $3$-cube.}
\label{fig:3cube1}
\end{figure}

Here we recall star unfolding, source unfolding, and farthest point map on the $3$-cube.
The (boundary of the) $3$-cube $M = \partial I^3$ has six square facets,
\begin{equation*}
\begin{aligned}
\sfL &= \set{0} \times I \times I , \quad \mbox{(left)} \\
\sfR &= \set{1} \times I \times I , \quad \mbox{(right)} \\
\sfF &= I \times \set{0} \times I , \quad \mbox{(front)} \\
\sfB &= I \times \set{1} \times I , \quad \mbox{(back)} \\
\sfD &= I \times I \times \set{0} , \quad \mbox{(down)} \\
\sfU &= I \times I \times \set{1} , \quad \mbox{(up)} ,
\end{aligned}
\end{equation*}
where we adopt a notation like Rubik's cube.
Let $\iota(x,y,z) = (1-x,1-y,1-z)$.

Let $p = (a,b,0) \in \sD$, and assume 
\[ 0 \le b \le a \le 1/2.
\]
Let $S$ be the union of the shortest paths from $p$ to the eight vertices of $M$
(as in Figure~\ref{fig:3cube1}(a)).
Then $M \setminus S$ can be unfolded to a polygon $M_2 \subset \aff(\sU) \simeq \real^2$
as in Figure~\ref{fig:3cube1}(b).

\begin{figure}
  \centering
 \includegraphics[width=40mm]{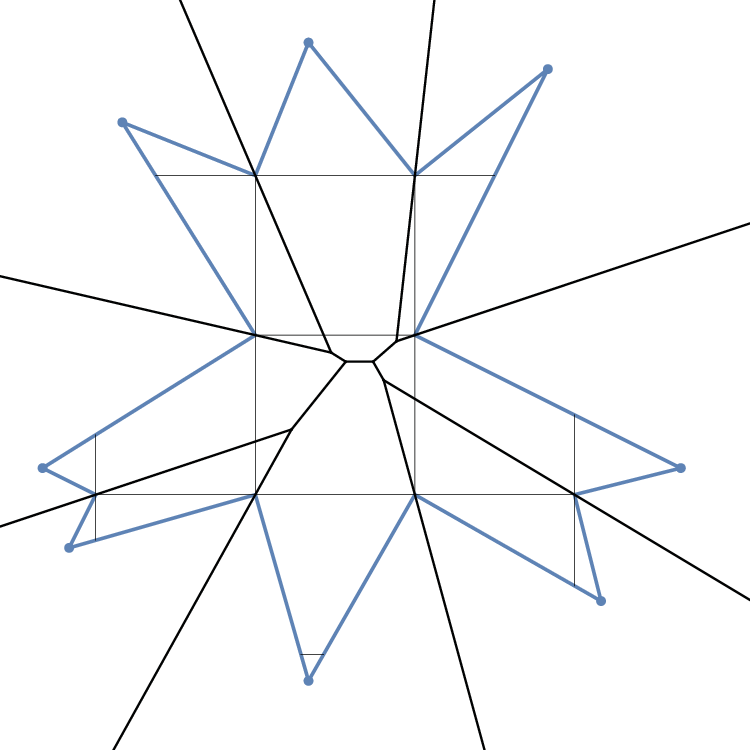}
    \label{fig:3cube-voronoi}
\includegraphics[width=40mm]{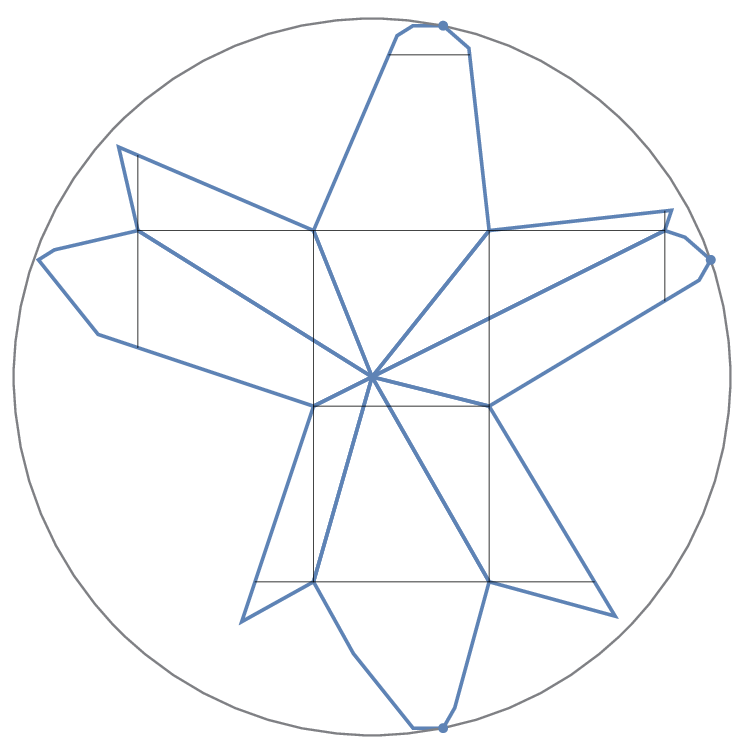}
    \label{fig:3cube-source}
  \label{fig:3cube2}
  \caption{(a) Voronoi tessellation of the star unfolding of the $3$-cube with respect to the eight source images.
  (b) A source unfolding of the $3$-cube.}
\end{figure}

The eight source images are defined by
\[ \src(\sU) = \set{ p_L \mid L \in \lbl(\sU)},
\quad 
\lbl(\sU) = \set{\sF, \sB, \sBL, \sL, \sLF, \sBR, \sR, \sRF},
\]
where $p_L = \phi_{\sU L \sD}(p)$.
They are given as
\begin{align*}
p_\sF &= (a, -1-b,1) , 
\quad 
p_\sB = (a, 3-b,1) , \\
p_\sBL &= (-1+b,2+a,1) , 
\quad 
p_\sL = (-1-a,b,1) , 
\quad
p_\sLF = (-1-b,-a,1) , \\
p_\sBR &= (2-b,3-a,1) ,
\quad
p_\sR = (3-a,b,1) , 
\quad
p_\sRF = (2+b, -1+a, 1) .
\end{align*}
In the intrinsic metric $\dist$ on $M$ we have 
\[ \dist(p,q) = \min_{p_L \in \src(\sU)} |q - p_L|
\]
for $q \in \sU$.
Let 
\[ V_F = \set{ w \in M_2 \mid |w - p_F| \le |w - p_{F'}|, 
\forall F' \in \lbl(\sU)}
\]
be the Voronoi region of the source image $p_F \in \src(\sU)$,
so $M_2 = \bigcup V_F$,
Figure~\ref{fig:3cube-voronoi}(a).
Let
\[ W_F := (\phi_{\sU F \sD})^{-1}(V_F) 
  \subset \aff(\sD),
\quad
F \in \lbl(\sU).
\]
Then $M_1 = \bigcup W_F$ is a source unfolding of $M$,
 Figure~\ref{fig:3cube-source}(b).

\begin{proposition}
The limit set of the farthest point map on $M = \partial I^3$ is the union of the diagonals of the six square facets.
\end{proposition}
\begin{proof}
Let $p=(a,b,0)$, $0 \le b \le a \le 1/2$,
and consider a sequence $\set{q_i}_i$ defined by
$q_0 = p$, 
$q_{i+1} \in \iota(f(q_i))$,
where $\iota (x,y,z) := (1-x, 1-y, 1-z)$.
First suppose that $a < 1/2$.
The farthest point of $p$ is 
the circumcenter of $p_\sF$, $p_\sR$, $p_\sB$
and is given by
\[ f(p) = \left(1 - \frac{a+2b(1-b)}{3-2a}, 1-b, 1 \right) ,
\]
where for simplicity a one-point set is identified with its element.
We have $\iota f(p) = \left( \frac{a+2b(1-b)}{3-2a}, b, 0 \right)$, so
\[ \lim_{k\to\infty} (\iota f)^k(p) = (b,b,0).
\]
If $0 \le b < a=1/2$, there are two farthest points
which are symmetric with respect to the plane $H_{a=1/2}=\set{(1/2,y,z)}$.
The above argument applies to both of them.
\end{proof}

\section{Star and source unfolding of the $4$-cube}

From this section on, we will deal with the (boundary of the) $4$-cube.
The $4$-cube $I^4$ has eight facets
\begin{align*}
\sfL &= \set{0} \times I \times I \times I, \quad \mbox{(left)} \\
\sfR &= \set{1} \times I \times I \times I, \quad \mbox{(right)} \\
\sfF &= I \times \set{0} \times I \times I, \quad \mbox{(front)} \\
\sfB &= I \times \set{1} \times I \times I, \quad \mbox{(back)} \\
\sfD &= I \times I \times \set{0} \times I, \quad \mbox{(down)} \\
\sfU &= I \times I \times \set{1} \times I, \quad \mbox{(up)} \\
\sfS &= I \times I \times I \times \set{0}, \quad \mbox{(source)} \\ 
\sfG &= I \times I \times I \times \set{1}, \quad \mbox{(goal)} 
\end{align*}
each of which is isometric to the $3$-cube.
Let $\iota(x,y,z,w)=(1-x,1-y,1-z,1-w)$ by abuse of notation.

Let
\[
 \Delta := \set{ (a,b,c,0) \mid 0 \le c \le b \le a \le 1/2}.
\]
Let $H_{a=\frac{1}{2}} = \set{(a,b,c,d) \mid a=1/2}$,
$H_{a=b} = \set{(a,b,c,d) \mid a=b}$,
$H_{b=c} = \set{(a,b,c,d) \mid b=c}$,
$H_{c=0} = \set{(a,b,c,d) \mid c=0} \subset \real^4$.
Lemma~\ref{3} implies that 
\begin{equation}
\label{5}
f(\Delta^\circ) \subset \iota \Delta.
\end{equation}
On the boundary $\partial \Delta$, 
it will be shown later that 
$\iota(f(\Delta \cap H)) \subset \Delta \cap H$
for $H=H_{a=b}$,
$H_{b=c}$, $H_{c=0}$.
For $H= H_{a=\half}$, 
it will be shown that
$\iota(f(\Delta \cap H)) \subset \Delta \cup r_H(\Delta)$.

Let $p = (a,b,c,0) \in \Delta$ be a source point.
Let
\begin{align*}
\src(\sfS) &= \set{p},
\\
\src(\sfD) &= \set{\phi_{\sfD\sfS}(p)},
\\
\src(\sfF) &= \set{\phi_{\sfF\sfS}(p), \phi_{\sfF\sfD\sfS}(p)},
\\
\src(\sfL) &= \set{\phi_{\sfL\sfS}(p), \phi_{\sfL\sfD\sfS}(p),
\phi_{\sfL\sfF\sfS}(p),\phi_{\sfL\sfF\sfD\sfS}(p)},
\\
\src(\sfR) &= \set{\phi_{\sfR\sfS}(p),\phi_{\sfR\sfD\sfS}(p),
\phi_{\sfR\sfF\sfS}(p), \phi_{\sfR\sfF\sfD\sfS}(p)},
\\
\src(\sfB) &= \set{\phi_{\sfB\sfS}(p),\phi_{\sfB\sfD\sfS}(p),
\phi_{\sfB\sfL\sfS}(p),\phi_{\sfB\sfL\sfD\sfS}(p),
\phi_{\sfB\sfR\sfS}(p),\phi_{\sfB\sfR\sfD\sfS}(p)},
\\
\src(\sfU) &= 
\set{ \phi_{L\sS}(p) \mid
L \in \set{\sU,\sUF,\sUB,\sUL,\sULF,\sUBL,\sUR,\sURF,\sUBR} },
\\
 \src(\sfG) &= \set{p_L \mid L \in \lbl(\sG)},
\end{align*}
where $p_L := \phi_{\sfG L \sfS}(p)$ and 
\begin{align*}
\lbl(\sG) &= \set{ \sU, \sD, 
\sUF, \sfF, \sfF\sfD, 
\sU\sfB, \sfB, \sfB\sfD,
\sU\sfL, \sfL, \sfL\sfD, 
\sU\sfR, \sfR, \sfR\sfD, 
\sU\sfL\sfF, \sfL\sfF, \sfL\sfF\sfD,
\\
& 
\qquad 
\sU\sfR\sfF, \sfR\sfF, \sfR\sfF\sfD,
\sU\sfB\sfL, \sfB\sfL, \sfB\sfL\sfD,
\sU\sfB\sfR, \sfB\sfR, \sfB\sfR\sfD
}.
\end{align*}
For each   
facet $F \in \set{\sfS, \sfG, \sfU,\sfD, \sfF, \sfB, \sfL, \sfR}$,
we have $\src(F) \subset \aff(F)$, 
and
\begin{equation}
\label{11}
\dist(p,q) 
= \min_{p' \in \src(F)} |p' - q| ,
\quad \forall q \in F .
\end{equation}
By Lemma~\ref{3}, we have $f(p) \subset \iota \Delta \subset \sG$.
The set $\src(\sfG)$ consists of
$26 = 3^3-1$ source points
\begin{align*}
p_\sfU &= (a,b,3-c,1), \quad
p_\sfD = (a,b,-1-c,1), \\
p_{\sUF} &= (a,-1+c,2+b ,1), \quad
p_\sF = (a,-1-b,c ,1), \quad
p_{\sFD} = (a,-1-c,-b ,1), \\
p_{\sUB} &= (a,2-c,3-b,1), \quad
p_\sB = (a,3-b,c,1), \quad
p_{\sBD} = (a,2+c,-1+b,1), \\
p_{\sUL} &= (-1+c,b,2+a ,1), \quad
p_\sL = (-1-a,b,c ,1), \quad
p_{\sLD} = (-1-c,b,-a ,1), \\
p_{\sUR} &= (2-c,b,3-a,1), \quad
p_\sR = (3-a,b,c ,1), \quad
p_{\sRD} = (2+c,b,-1+a ,1), \\
p_{\sULF} &= (-1+c,-a,2+b ,1), \quad
p_{\sLF} = (-1-b,-a,c ,1), \quad
p_{\sLFD} = (-1-c,-a,-b ,1), \\
p_{\sURF} &= (2-c,-1+a,2+b ,1), \quad
p_{\sRF} = (2+b,-1+a,c ,1), \\
& \qquad p_{\sRFD} = (2+c,-1+a,-b ,1), \\
p_{\sUBL} &= (-1+b,2-c,2+a ,1), \quad
p_{\sBL} = (-1+b, 2+a, c ,1), \\
& \qquad p_{\sBLD} = (-1+b, 2+c, -a ,1), \\
p_{\sUBR} &= (2-b,2-c,3-a ,1) , \quad
p_{\sBR} = (2-b,3-a,c ,1), \\
& \qquad p_{\sBRD} = (2-b,2+c, -1+a ,1)
.
\end{align*}
We will often identify $\aff(\sG)$ with $\real^3$ 
by identifying $(x,y,z,1)$ with $(x,y,z)$.

\begin{figure}
  \centering
\includegraphics[width=52mm]{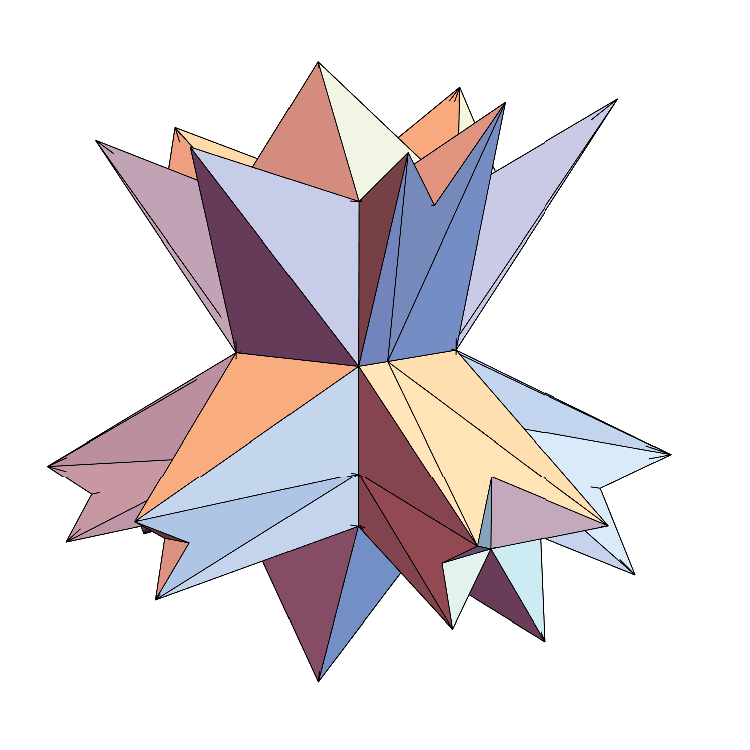}
    \label{fig:4cube-star}
\includegraphics[width=52mm]{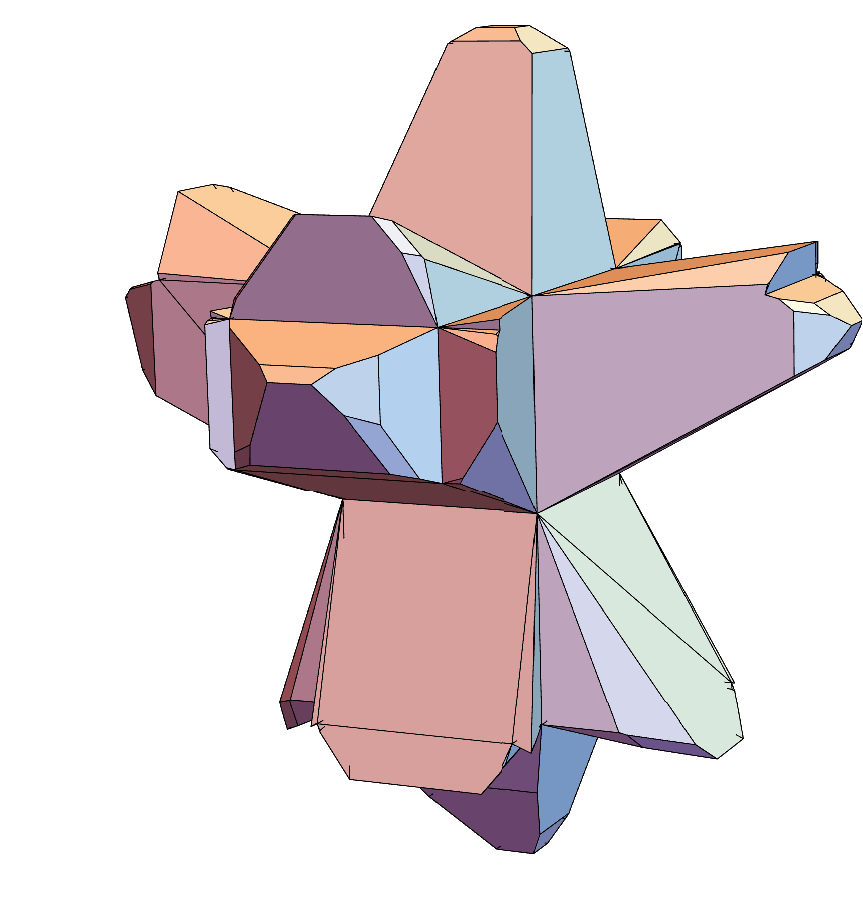}
    \label{fig:4cube-source}
  \label{fig:8OSMElogosubfigs}
  \caption{(a) A star unfolding of the $4$-cube
  with $26=3^3-1$ source images of $p=(2/5,1/3,1/6,0)$.
 The concave quadrilateral faces are additional cuts
 whose diagonals are cut points.
 (b) A source unfolding of the $4$-cube.}
\end{figure}

There exists a star unfolding
$M \setminus S \to M_2 \subset \aff(\sG) \simeq \real^{3}$ for some subset $S \subset M$, 
 with the source image $\src(\sfG)$ as vertices,
  Figure~\ref{fig:4cube-star}.
In this figure, 
each triangular face is a union of shortest paths from a source point to warped points. 
The concave quadrilateral faces are additional cuts 
that are not given in the classical definition of the star unfolding.
The diagonal of a concave quadrilateral face
consists of cut points.

Let 
\[ V_L = \set{ w \in M_2 \mid |w - p_L| \le |w - p_{L'}|,
 \forall L' \in \lbl(\sG)}
\]
be the Voronoi domain of the source image $p_L \in \src(\sG)$.
We have $M_2 = \bigcup V_L$.
Let
\[ W_L := (\phi_{\sG L \sS})^{-1}(V_L) 
 = \phi_{\sS L^{-1} \sG} (V_L) \subset \aff(\sS).
\]
Then $M_1 = \bigcup W_L$ is a source unfolding of $M$,
see Figure~\ref{fig:4cube-source}.

\section{Delaunay triangulations in the planes $H_{x=a}$, $H_{y=b}$}

In Sections 4-6,
we consider the Voronoi tessellation of $\aff(\sG)$
with respect to $\src(\sG)$.

A farthest point is a Voronoi point, which is the circumcenter of 
a Delaunay polyhedron.  
We will show that any farthest point
is the circumcenter of a polyhedron 
whose vertices belong to
\begin{equation}
\label{13}
 \set{p_F \mid F \in \BackRight} \subset \src(\sG), 
\end{equation}
where
\[ \BackRight = \set{\sU, \sD, \sUR, \sR, \sRD, \sUB, \sB, \sBD,
 \sUBR, \sBR, \sBRD} .
\]

In preparation, we consider Delaunay triangulations of 
the subsets
\begin{align*}
\src(\sfG) \cap H_{x=a} 
&:= \set{ \pU, \pD, \pUF, \pF, \pFD, \pUB, \pB, \pBD} ,
\\
\src(\sfG) \cap H_{y=b} 
& :=  \set{\pU, \pD, \pUL, \pL, \pLD, \pUR, \pR, \pRD} ,
\end{align*}
in the planes
$H_{x=a} = \set{(x,y,z,w) \mid x=a}$ and 
$H_{y=b} = \set{(x,y,z,w) \mid y=b}$.
We will show that their Delaunay triangles 
in the planes $H_{x=a}$, $H_{y=b}$ 
are also Delaunay triangles in the Delaunay tetrahedralization
of $\src(\sG)$
in $\aff(\sG)$.
By Lemma~\ref{3}, 
we can assume that the farthest point(s) 
are in the region $\iota\Delta$.
The above Delaunay triangles act like walls,
so it is sufficient to consider the site points in (\ref{13}) 
as the vertices of the Delaunay polyhedron for the farthest point.

First consider the plane $H_{y=b}$.
Let $\pi_{xz} : \real^4 \to \real^2$,
$\pi_{xz}(\xi,\eta,\zeta,\omega) = (\xi,\zeta)$, be a projection.
Let
\[
 W_{y=b} 
:= \pi_{xz}( \src(\sfG) \cap H_{y=b})
= 
\set{\pi_{xz}(p_F) \mid 
F\in\set{\sfU, \sfD,\sfU\sfL, \sfL,\sfL\sfD, \sfU\sfR,\sfR,\sfR\sfD}
} .
\]
Denote by $q_F := \pi_{xz}(p_F)$.
We will give the Delaunay triangulation of $W_{y=b}$ in $\real^2$.
Let $\Delta_{0}=\set{(a,c) \mid 0 \le c \le a \le 1/2}$.

\begin{lemma}
The line segment $[\qU,\qD]$ is a Delaunay edge of $W_{y=b}$
in $\real^2$. 
\end{lemma}
\begin{proof}
We show that the circle $C$ with the center $(1-a,1-c)$ passing through
$\qU=(a,3-c)$, $\qD=(a,-1-c)$ is an empty circle.
That is, the other points 
$\qUL, \qL, \qLD, \qUR, \qR, \qRD \in W_{y=b}$ are outside or on
$C$. 

Let
\begin{align*}
\varphi_\sUL(a,c) &= \frac{(1+c)(1-a-c)}{1+a-c},
\\
\varphi_\sL(a,c) &= -\frac{1}{2}+\frac{(3-2c)(1+2c)}{2+4a},
\\
\varphi_\sLD(a,c) &= \frac{(1-a+c)(1-c)}{1+a+c},
\end{align*}
and $\varphi_\sUR(a,c) = 1 - \varphi_\sUL(1-a,c)$,
$\varphi_\sR(a,c) = 1 - \varphi_\sL(1-a,c)$,
$\varphi_\sRD(a,c) = 1 - \varphi_\sLD(1-a,c)$.
Then 
$(\varphi_F(a,c),1-c)$ is the circumcenter of $\triangle(\qU,\qD,q_F)$
for $F\in \set{\sUL, \sL, \sLD, \sUR,\sR,\sRD}$.

We have 
\[ \varphi_\sUL, \varphi_\sL, \varphi_\sLD
\le  1-a 
\le \varphi_\sUR, \varphi_\sR, \varphi_\sRD
\]
for $(a,c) \in \Delta_{0}$.
This shows that $C$ is an empty circle. 
\end{proof}

\begin{figure}
(a)
\includegraphics[scale=0.4]{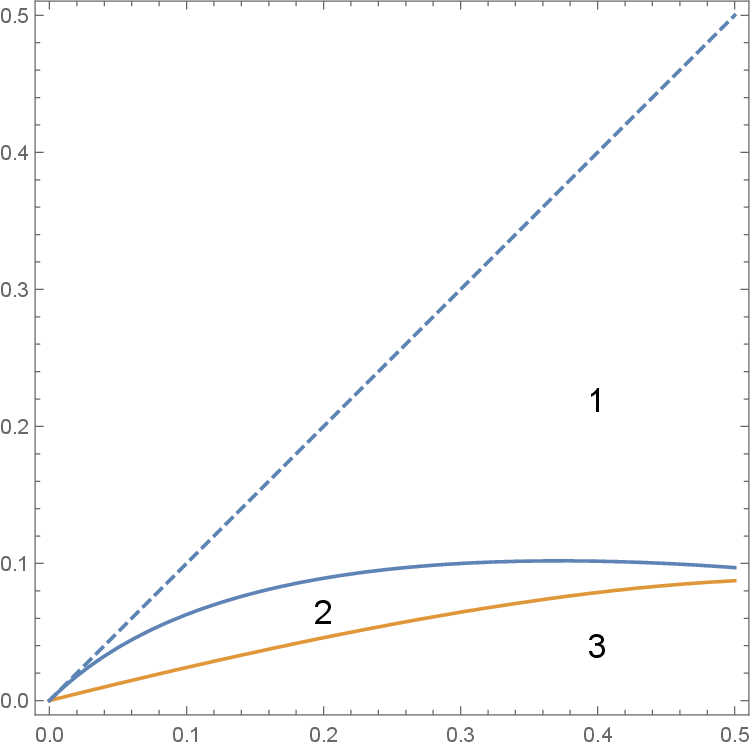}
(b)
\includegraphics[scale=0.4]{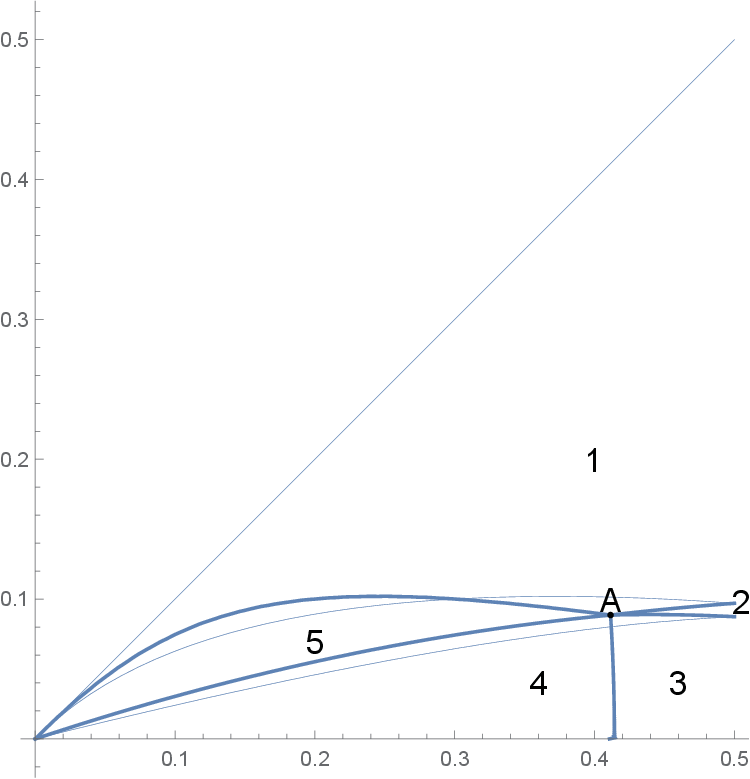}
\label{12}
\caption{(a) The regions $\Delta_i$, $i=1,2,3$.
(b) The regions $\Delta_{iL}$, $i=1,\dots,5$.}
\label{fig:bif}
\end{figure}

The Delaunay triangulation on the right side of $[\qU,\qD]$ is given as follows.
Let
\begin{align*}
 \psi_1(a,c) &:= \outcircle(\qU,\qD,\qR,\qRD) 
\\
  &= 8(-a+c + a^2 + 4 a c + c^2 - 2a^2c-2c^3), \\
 \psi_2(a,c) &:= \outcircle(\qU,\qRD,\qR, \qUR)
  \\
 &= 4(-2a + 8c + 3a^2 - 8ac - 5c^2 - a^3 + 3a^2c - ac^2 + 3c^3 ).
\end{align*}
Let
\begin{align*}
\Delta_1 &= \set{(a,c) \in \Delta_0 \mid \psi_1(a,c)\ge0}
\\
\Delta_2 &= \set{(a,c) \in \Delta_0 \mid \psi_1(a,c) \le 0 \le \psi_2(a,c)}
\\
\Delta_3 &= \set{(a,c) \in \Delta_0 \mid \psi_2(a,c) \le 0}.
\end{align*}
The regions $\Delta_i$, $i=1,2,3$, are shown in Figure~\ref{12}(a).
If $(a,c) \in \Delta_1$, then
$\triangle(\qU,\qD,\qR)$,
$\triangle(\qD,\qR,\qRD)$,
$\triangle(\qU,\qUR,\qR)$ are Delaunay triangles.
If $(a,c) \in \Delta_2$,
then 
$\triangle(\qU,\qUR, \qR)$,
$\triangle(\qU,\qR,\qRD)$,
$\triangle(\qU,\qD,\qRD)$ are Delaunay triangles.
If $(a,c) \in \Delta_3$,
then
$\triangle(\qU,\qD,\qRD)$,
$\triangle(\qU,\qUR,\qRD)$,
$\triangle(\qUR,\qR,\qRD)$ are Delaunay triangles.

The Delaunay triangulation on the left side of $[\qU,\qD]$ is given as follows (which will no be used later in this paper).
Let
\begin{align*}
\psi_{1L}(a,c) &:= \outcircle(\qL,\qLD,\qD,\qU)
 = \psi_1(1-a,c)
\\
\psi_{2L}(a,c) &:= \outcircle(\qU,\qUL,\qL,\qLD)
=  \psi_2(1-a,c)
\\
\psi_{3L}(a,c) &:= \outcircle(\qLD,\qD,\qU,\qUL)
 = 16c(-1+2a+a^2+c^2)
\\
\psi_{4L}(a,c) &:= \outcircle(\qUL,\qL,\qLD,\qD)
= 4(a - a^3 - 3c - 2ac+a^2c - 2c^2 - ac^2 + c^3)
\\
\psi_{5L}(a,c) &:= \outcircle(\qD,\qU,\qUL,\qL)
= 8 (a - a^2 - c - 4 a c - 2a^2c +3c^2 -2c^3).
\end{align*}

Let $A = (-1+\sqrt{7})/4, (3-\sqrt{7})/4)$.
Let
\begin{align*}
\Delta_{1L} &= \set{(a,c) \in \Delta_0 \mid c \ge (3+\sqrt{7})/4, 
\psi_{5L}(a,c) \le 0 \le \psi_{1L}(a,c)}
\\
\Delta_{2L} &= \set{(a,c) \in \Delta_0 \mid a \ge (-1+\sqrt{7})/4,
\psi_{1L}(a,c) \le 0 \le \psi_{2L}(a,c)}
\\
\Delta_{3L} &= \set{(a,c) \in \Delta_0 \mid a \ge (-1+\sqrt{7})/4,
\psi_{2L}(a,c) \le 0 \le \psi_{3L}(a,c)}
\\
\Delta_{4L} &= \set{(a,c) \in \Delta_0 \mid c \le (3-\sqrt{7})/4,
\psi_{3L}(a,c) \le 0 \le \psi_{4L}(a,c)}
\\
\Delta_{5L} &= \set{(a,c) \in \Delta_0 \mid a \le (-1+\sqrt{7})/4,
\psi_{4L}(a,c) \le 0 \le \psi_{5L}(a,c)}.
\end{align*}
These regions are shown in Figure~\ref{fig:bif}(b).
If $(a,c) \in \Delta_{1L}$,
then $\triangle(\qU,\qD,\qL)$,
$\triangle(\qD,\qL,\qLD)$,
$\triangle(\qU,\qUL,\qL)$ are Delaunay triangles.
If $(a,c) \in \Delta_{2L}$,
then $\triangle(\qU,\qUL,\qL)$,
$\triangle(\qU,\qL,\qLD)$,
$\triangle(\qU,\qD,\qLD)$ are Delaunay triangles.
If $(a,c) \in \Delta_{3L}$,
then $\triangle(\qU,\qD,\qLD)$,
$\triangle(\qU,\qUL,\qLD)$,
$\triangle(\qUL,\qL,\qLD)$ are Delaunay triangles.
If $(a,c) \in \Delta_{4L}$,
then $\triangle(\qUL,\qL,\qLD)$,
$\triangle(\qD,\qUL,\qLD)$,
$\triangle(\qU,\qD,\qUL)$ are Delaunay triangles.
If $(a,c) \in \Delta_{5L}$,
then $\triangle(\qU,\qD,\qUL)$,
$\triangle(\qD,\qUL,\qL)$,
$\triangle(\qD,\qL,\qLD)$ are Delaunay triangles.

Let $\pi_{yz}(\xi,\eta,\zeta,\omega) = (\eta,\zeta)$.
The Delaunay triangulation of the set
\[
 W_{x=a} := \set{\pi_{yz}(p_F) \mid 
F \in \set{\sU, \sD, \sUF, \sF, \sFD, \sUB, \sB, \sBD } }
\subset \real^2
\]
is similar to that of $W_{y=b}$.
Denote by $q'_F = \pi_{yz}(p_F)$.
Let $\Delta'_0 = \set{(b,c) \mid 0 \le c \le b \le 1/2} \simeq \Delta_0$.
Let
$\Delta'_1 = \set{(b,c) \in \Delta'_0 \mid \psi_1(b,c)\ge0}$,
$\Delta'_2 = \set{(b,c) \in \Delta'_0 \mid \psi_1(b,c) \le 0 \le \psi_2(b,c)}$,
$\Delta'_3 = \set{(b,c) \in \Delta'_0 \mid \psi_2(b,c) \le 0}$.
The line segment $[\qU',\qD']$ is a Delaunay edge 
of $W_{x=a}$.
The Delaunay triangulation of $W_{x=a}$ in $H_{x=a}$ has the following
Delaunay triangles (in the back side of $[\qU',\qD']$).
If $(b,c) \in \Delta_1$, then
$\triangle(\qU',\qD',\qB')$,
$\triangle(\qD',\qB',\qBD')$,
$\triangle(\qU',\qUB',\qB')$ are Delaunay triangles.
If $(b,c) \in \Delta_2$,
then 
$\triangle(\qU',\qUB', \qB')$,
$\triangle(\qU',\qB',\qBD')$,
$\triangle(\qU',\qD',\qBD')$ are Delaunay triangles.
If $(b,c) \in \Delta_3$,
then
$\triangle(\qU',\qD',\qBD')$,
$\triangle(\qU',\qUB',\qBD')$,
$\triangle(\qUB',\qB',\qBD')$ are Delaunay triangles.

\section{Delaunay triangles in $H_{x=a}$, $H_{y=b}$}

In this section we show that some Delaunay triangles
in the Delaunay triangulations of $\src(\sG) \cap H_{x=a}$, $\src(\sG) \cap H_{y=b}$
in $H_{x=a}$, $H_{y=b}$,
are also Delaunay triangles
in the Delaunay tetrahedralization of $\src(\sG)$ in $\aff(\sG)$.

Let
\begin{align*}
\Delta_{1*} &= \set{(a,b,c,0) \in \Delta  \mid \psird(a,c)\ge 0 }, \\
\Delta_{2*} &= \set{(a,b,c,0) \in \Delta  \mid \psird(a,c) \le 0 \le \psiur(a,c)},
\\
\Delta_{3*} &= \set{(a,b,c,0) \in \Delta  \mid \psiur(a,c) \le 0}
 \\
\Delta_{*1} &= \set{(a,b,c,0) \in \Delta  \mid \psird(b,c) \ge 0 }, \\
\Delta_{*2} &= \set{(a,b,c,0) \in \Delta  \mid \psird(b,c) \le 0 \le \psiur(b,c)},
\\
\Delta_{*3} &= \set{(a,b,c,0) \in \Delta  \mid \psiur(b,c) \le 0},
\end{align*}
and
$\Delta_{ij} := \Delta_{i*} \cap \Delta_{*j}$,
$ij \in \set{11,12,21,22,31,32,33}$.
Note that $\Delta_{1*} \cap \Delta_{*3}$,
$\Delta_{2*} \cap \Delta_{*3}$ are empty sets.

Let $q(F_1,F_2,F_3)$ be the circumcenter of 
$\triangle(p_{F_1}, p_{F_2}, p_{F_3})$
for $(F_1,F_2,F_3) = (\sU,\sD,\sR)$,
$(\sD,\sR,\sRD)$,
$(\sU,\sUR,\sR)$,
$(\sU,\sR,\sRD)$,
$(\sU,\sD,\sRD)$,
$(\sU,\sUR,\sRD)$,
$(\sUR,\sR,\sRD)$.
So
\begin{align*}
q(\sU,\sD,\sR) &=
\left( \frac{3}{2} -\frac{(3-2c)(1+2c)}{2(3-2a)}, b, 1-c, 1 \right)
\\
q(\sD,\sR,\sRD)
 &=
 \left( \frac{3}{2} -\frac{(1-a)(1+2c)}{2(1-2a+a^2+c+c^2)}, b, 
  -\frac{1}{2} + \frac{(1-a)(3-2a)}{2(1-2a+a^2+c+c^2)} , 1 \right)
\\
q(\sU,\sUR,\sR)
&= \left(
\frac{3}{2}+\frac{(1-c)(3-2c)}{2 (3 - 3 a + a^2 - 2 c + c^2)}, b,
\frac{3}{2}-\frac{(1-c)(3-2a)}{2 (3 - 3 a + a^2 - 2 c + c^2)}, 1 \right)
\\
q(\sU,\sR,\sRD)
&=
\left(
\frac{3}{2} - \frac{(1-a+2c)(3-2c)}{2 (3 - 4 a + a^2 - c + c^2)}, b,
\frac{3}{2} - \frac{(1-a+2c)(3-2a)}{2 (3 - 4a + a^2 - c + c^2)}, 1 \right)
\\ 
q(\sU,\sD,\sRD)
 &=
 \left(1 -\frac{(1-c)(a+c)}{2-a+c}, b, 1-c, 1 \right)
\\
q(\sU,\sUR,\sRD)
&=
\left(
 2-\frac{2(2-a)(1-c)}{4-4a+a^2-2c+c^2}, b,
 1- \frac{2c(1-c)}{4-4a+a^2-2c+c^2}, 1 \right)
\\ 
q(\sUR,\sR,\sRD)
&=
\left( 2 -\frac{(2-a)(1-a+c)}{2-3a+a^2-c+c^2}, b,
 1- \frac{c(1-a+c)}{2-3a+a^2-c+c^2}, 1 \right) .
\end{align*}

\begin{lemma}
\label{4}
If $(a,b,c) \in \Delta_{1*}$, then
$\triangle(\pU,\pD,\pR)$,
$\triangle(\pD,\pR,\pRD)$,
$\triangle(\pU,\pUR,\pR)$ are Delaunay triangles 
in the Delaunay tetrahedralization of $\src(\sG)$ in $\aff(\sG)$.
If $(a,b,c) \in \Delta_{2*}$,
then 
$\triangle(\pU,\pUR, \pR)$,
$\triangle(\pU,\pR,\pRD)$,
$\triangle(\pU,\pD,\pRD)$ are Delaunay triangles.
If $(a,b,c) \in \Delta_{3*}$,
then
$\triangle(\pU,\pD,\pRD)$,
$\triangle(\pU,\pUR,\pRD)$,
$\triangle(\pUR,\pR,\pRD)$ are Delaunay triangles.
\end{lemma}

\begin{proof}
Let $\alpha_{b}(\xi,b,\zeta,1) = (\xi, 1-b, \zeta,1)$,
$\alpha_{b,c}(\xi,b,\zeta,1) = (\xi, c-b+\zeta, \zeta,1)$. 

First we show that the sphere $C$ with the center 
$\alpha = \alpha_b(q(\sU,\sD,\sR))$
passing through the points $\pU, \pD,\pR$ is an empty sphere
if $\psi_1(a,c)\ge0$.
That is, the other points of $\src(\sG)$ are outside or on the boundary of $C$.
In fact, 
denote the circumcenter of $\pU, \pD, \pR, p_F$
by
\[ \left( \frac{3}{2} -\frac{(3-2c)(1+2c)}{2(3-2a)}, \phi_F, 1-c , 1 \right),
\] 
for $F$.  Then we have
\[ \phi_F \le 1-b
\]
for $F \in \Front := \set{\sUF, \sF,\sFD, \sULF, \sLF, \sLFD, \sURF, \sRF, \sRFD}$,
and
$1-b \le \phi_F$ for 
$F \in \Back := \set{ \sUB, \sB, \sBD, \sUBL, \sBL, \sBLD, \sUBR, \sBR, \sBRD}$.

Similarly, we can show the followings.
If $\psi_1(a,c) \ge 0$,
the sphere with the center 
$\alpha = \alpha_b(q(\sD,\sR,\sRD))$
passing through $\pD, \pR, \pRD$ is an empty sphere. 

If $\psi_2(a,c)\ge 0$, 
the sphere with the center
$\alpha_{b,c}(q(\sU,\sUR,\sR))$
passing through $\pU, \pUR, \pR$ is an empty sphere.

If $\psi_1(a,c) \le 0 \le \psi_2(a,c)$,
the sphere with center $\alpha_{b}(q(\sU,\sR,\sRD))$
passing through $\pU,\pRD,\pR$ is an empty sphere.

If $\psi_1(a,c) \le 0$,
the sphere with center 
$\alpha_b(q(\sU,\sD,\sRD))$
passing through $\pU,\pD,\pRD$ is an empty sphere.

If $\psi_2(a,c) \le 0$, 
the sphere with center 
$\alpha_{b,c}(q(\sU,\sUR,\sRD))$
passing through $\pU,\pRD,\pUR$ is an empty sphere,
and the sphere  
$\alpha_{b,c}(q(\sUR,\sR,\sRD))$
passing through $\pRD,\pR,\pUR$ is an empty sphere.
\end{proof}

\begin{lemma}
\label{16}
If $(a,b,c) \in \Delta_{*1}$, then
$\triangle(\pU,\pD,\pB)$,
$\triangle(\pD,\pB,\pBD)$,
$\triangle(\pU,\pUB,\pB)$ are Delaunay triangles
in the Delaunay tetrahedralization of $\src(\sG)$ in $\aff(\sG)$.
If $(a,b,c) \in \Delta_{*2}$,
then 
$\triangle(\pU,\pUB, \pB)$,
$\triangle(\pU,\pB,\pBD)$,
$\triangle(\pU,\pD,\pBD)$ are Delaunay triangles.
If $(a,b,c) \in \Delta_{*3}$,
then
$\triangle(\pU,\pD,\pBD)$,
$\triangle(\pU,\pUB,\pBD)$,
$\triangle(\pUB,\pB,\pBD)$ are Delaunay triangles.
\end{lemma}
\begin{proof}
The proof is similar to that of Lemma~\ref{4}.
\end{proof}

Now we can show the following Proposition.

\begin{proposition}
\label{15}
Let $q\in \iota\Delta$ be a farthest point of $p \in \Delta$.
Then $q$ is the circumcenter of a polyhedron whose vertices belong
to
\[ \set{ p_F \mid F \in \BackRight} \subset \src(\sG).
\]
\end{proposition}

To prove Proposition~\ref{15}, we prepare the following lemma.

\begin{lemma}
\label{14}
Suppose that a farthest point $q$ of $p$ is
the circumcenter of a Delaunay polyhedron $K$,
and suppose that $K$ is maximal in the sense that
$\partial C \cap K = \partial C \cap \src(\sG)$
for the circumsphere $C$ of $K$.
Then $q$ is inside $K$.
\end{lemma}
\begin{proof}
If $q$ is outside or on
the boundary of $K$, then there exists a plane $H$ passing through $q$ 
such that $K$ is contained in a (closed) half-space with respect to $H$. 
Let $v$ be a normal vector of $H$.  Then, by moving from $q$ in the 
direction of $v$, we obtain a point farther than $q$.
\end{proof}

\begin{proof}[Proof of Proposition~\ref{15}]
By Lemma~\ref{14}, 
$q$ is the circumcenter of a 
polyhedron $K$ whose
vertices belong to $\src(\sG)$.
With Lemmas~\ref{4}, \ref{16},
and $q \in \iota\Delta \subset \set{(x,y,z,1) \mid x\ge a, y\ge b}$,
the vertices of $K$
can be chosen as belonging to 
$\set{p_F \mid F \in \BackRight}$.
\end{proof}

\section{Corners of the Voronoi domain of $\pU$}

In this section we show that any farthest point $q \in \iota\Delta$
is a corner of the Voronoi domain of $\pU$.

Denote the tetrahedron 
with the vertices  $\pU, F_1,F_2,F_3 \in \BackRight$
by $t(F_1,F_2,F_3)$ where $\pU$ is omitted.
Denote the circumcenter of $t(F_1,F_2,F_3)$ by $c(F_1,F_2,F_3)$. 
We have
\begin{align*}
c(\sB, \sD, \sR) &=
 \left( 1 - \varphi_{1,c}(a), 1 - \varphi_{1,c}(b), 1-c
 \right),
 \quad \varphi_{1,c}(\xi) = \frac{\xi+2c(1-c)}{3-2\xi}
 \\
 c(\sB,\sD,\sRD) &=
 \left(
  1 - \varphi_{2,c}(a), 1 - \varphi_{1,c}(b), 1-c
 \right)
 \\
c(\sBD, \sD,\sR) &=
 \left(
  1 - \varphi_{1,c}(a), 1-\varphi_{2,c}(b), 1-c
 \right),
 \quad \varphi_{2,c}(\xi) = \frac{(\xi-c)(1+c)}{2-\xi-c}
 \\
c(\sBD,\sD,\sRD) &=
 \left(
  1 - \varphi_{2,c}(a), 1 - \varphi_{2,c}(b), 1-c
 \right)
 \\
\end{align*}

\begin{align*}
c(\sB,\sRD,\sR)
&= \bigg( \frac{3}{2} - \frac{(1-a+2c)(3-2c)}{2 (3 - 4 a -c + a^2 + c^2)}, 
\frac{3}{2} - \frac{(3-2a)(1-a+2c)(3-2c)}{2 (3 - 2 b) (3 - 4 a + a^2 - c + c^2)},
\\
& \qquad
\frac{3}{2} - \frac{(3-2a)(1-a+2c)}{2 (3 - 4 a + a^2 - c + c^2)}
\bigg) ,
\end{align*}
\begin{align*}
c(\sB,\sfRD,\sfUR)
 &= \bigg(1 - \frac{(2 - a + c) (a - c)}{4 - 4 a + a^2 - 2 c + c^2}, 
 \frac{3}{2} - \frac{(3-2c)(4-4a+2c+a^2-3c^2)}{(3 - 2 b) (4 - 4 a + a^2 - 2 c + c^2)},
\\
& \quad
  1 - \frac{2 (1 - c) c}{4 - 4 a + a^2 - 2 c + c^2}
\bigg) ,
\end{align*}
\begin{align*}
c(\sB,\sfBD,\sRD)
& = 
\bigg( 
3 - \frac{(4 - a - c) (3 - 3 b + c) (2 - b + c)}{2 (2 - a + c) (3 - 4 b + b^2 - c + c^2)},
 2 - \frac{(1 - b + c) (6 - b - 3 c) }{2 (3 - 4 b + b^2 - c + c^2) }, 
\\
& \quad 
\frac{(1 - b - c) (6 - b - 3 c) }{2 (3 -  4 b + b^2 - c + c^2)} 
\bigg), 
\end{align*}
\begin{align*}
c(\sB,\sfBD,\sUR)
& = 
\bigg(1 - \frac{(a - c) (3 - 3 b + c) (2 - b + c)}
{2 (2 - a - c) (3 - 4 b + b^2 - c + c^2)}, 
 \frac{3}{2} - \frac{(1-b+2c)(3-2c)}{2 (3 - 4 b - c + b^2 + c^2)},
\\
& \quad
 \frac{3}{2}- \frac{(3-2b)(1-b+2c)}{2 (3 - 4 b - c + b^2 + c^2)}
\bigg) ,
\end{align*}
\begin{align*}
c(\sBD,\sfRD,\sfUR)
&= 
 \bigg( 1 - \frac{(2 - a + c) (a - c)}{4 - 4 a + a^2 - 2 c + c^2}, 
 3 - \frac{(2-a+c)(2-a-c)(4-b-c)}{(2 - b + c) (4 - 4 a + a^2 - 2 c + c^2)},
\\
& \quad 
 1 - \frac{2 (1 - c) c}{4 - 4 a + a^2 - 2 c + c^2}
\bigg) ,
\end{align*}
\begin{align*}
c(\sfUB,\sfBD,\sfUR)
&= \bigg( 1 - \frac{(a - c) (2 - b + c) (2 - b - c)}
{(2 - a - c) (4 - 4 b + b^2 - 2 c + c^2)},
1 - \frac{(2 - b + c) (b - c)} {4 - 4 b + b^2 - 2 c + c^2},
\\
& \quad 
1 - \frac{2 (1 - c) c}{4 - 4 b + b^2 - 2 c + c^2} 
\bigg).
\end{align*}

Let
\begin{align*}
 \Delta_{22A} &= \set{p \in \Delta_{22} \mid 
\psi_{22}(a,b,c) \ge 0
} ,
\\
 \Delta_{22B} &= \set{p \in \Delta_{22} \mid
\psi_{22}(a,b,c) \le 0 } ,
\\
\Delta_{32A} &= \set{p\in \Delta_{32} \mid
\psi_{32}(a,b,c) \ge 0 } ,
\\
\Delta_{32B} &= \set{p \in \Delta_{32} \mid
\psi_{32}(a,b,c) \le 0 } ,
\end{align*}
where
\begin{align*}
\psi_{22}(a,b,c)
&:= \outsphere(\pU,\pB,\pBD,\pRD,\pR) / 4(a-b) 
\\
 & = 
 3 - 3 a - 3 b - 17 c + 3 a b + 8 a c + 8 b c + c^2
 - 4 a b c 
 - 2 a c^2 - 2 b c^2 + 4 c^3 
\\
&= - \frac{3-3b-8c+4bc+2c^2}{(a-b)(1-2c)} \psird(a,c)
   + \frac{3-3a-8c+4ac+2c^2}{(a-b)(1-2c)} \psird(b,c) ,
\\
\psi_{32}(a,b,c)
& := \outsphere(\pU,\pB,\pBD,\pRD,\pUR) \\
& = - \frac{4(a-b)(4-4a-6c+a^2+5c^2)}{3-3a-8c+4ac+2c^2}\psi_{22}(a,b,c)
\\
& \quad 
 + \frac{(2-a+c)(5c(1-c)(3-2c)(1+2c)-(3-4c)\psi_{22}(a,b,c))}{(3-3a-8c+4ac+2c^2)^2} \psi_2(a,c)
\end{align*}

Let
\begin{align*}
v_{11}(\sU) &=  \set{\rmc(\sB,\sD,\sR), \rmc(\sB,\sR,\sUR),  
\rmc(\sUB,\sB,\sUR)} ,
\\
v_{12}(\sU) &= \set{ \rmc(\sBD,\sD,\sR), \rmc(\sB,\sBD,\sR),  
\rmc(\sB,\sR,\sUR), \rmc(\sUB,\sB,\sUR) } ,
\\
v_{21}(\sU) &= \set{ \rmc(\sB,\sD,\sRD) ,
 \rmc(\sB,\sRD,\sR), 
 \rmc(\sB,\sR,\sUR),
 \rmc(\sUB,\sB,\sUR) } ,
\\
v_{22A}(\sU) &=  \set{ \rmc(\sBD,\sD,\sRD),
 \rmc(\sB,\sBD,\sRD), \rmc(\sB,\sRD,\sR),
 \rmc(\sB,\sR,\sUR),
 \rmc(\sUB,\sB,\sUR)} ,
\\
v_{22B}(\sU) &= \set{ \rmc(\sBD,\sD,\sRD),
 \rmc(\sBD,\sRD,\sR), \rmc(\sB,\sBD,\sR),
 \rmc(\sB,\sR,\sUR),
 \rmc(\sUB,\sB,\sUR)} ,
\\
v_{31}(\sU) &= \set{ \rmc(\sB,\sD,\sRD), 
\rmc(\sB,\sRD,\sUR), 
\rmc(\sUB,\sB,\sUR)} ,
\\
v_{32A}(\sU) &=  \set{\rmc(\sBD,\sD,\sRD), 
 \rmc(\sB,\sBD,\sRD), \rmc(\sB,\sRD,\sUR),
\rmc(\sUB,\sB,\sUR)} ,
\\
v_{32B}(\sU) &=  \set{\rmc(\sBD,\sD,\sRD), 
\rmc(\sBD,\sRD,\sUR), \rmc(\sB,\sBD,\sUR),
\rmc(\sUB,\sB,\sUR)} ,
\\
v_{33}(\sU) &=  \set{\rmc(\sBD,\sD,\sRD), 
\rmc(\sBD,\sRD,\sUR), 
\rmc(\sUB,\sBD,\sUR)} .
\end{align*}
Note that the above definitions coincide on the boundary,
e.g., $v_{11}(\sU) = v_{12}(\sU)$ if $p \in \Delta_{11} \cap \Delta_{12}$.

\begin{lemma}
Let $k \in \set{11,12,21,22A,22B,31,32A,32B,33}$.
If $p \in \Delta_k$, then each $q \in v_k(\sU)$ 
is a corner of the Voronoi domain of $\pU$ 
in the Voronoi tessellation of $\aff(\sG)$ with respect to
$\set{p_F \mid F \in \BackRight}$.
\end{lemma}

\begin{proof}
We identify $\aff(\sG)$ with $\real^3$.
 Let
\begin{align*}
f_{\sD}(x,y) &= 1 - c ,
\\
f_{\sUB}(x,y) &= (- 2 + 3 b - c + (2 - b - c) y)/(b - c) ,
\\
f_{\sB}(x,y) &= (3 b - 3 c + (3 - 2 b) y)/(3 - 2 c) ,
\\
f_{\sBD}(x,y) &= (2 + b - 5 c + (2 - b + c) y)/(4 - b - c) ,
\\
f_{\sUR}(x,y) &= (- 2 + 3 a - c + (2 - a - c) x)/(a - c) ,
\\
f_{\sR}(x,y) &= (3 a - 3 c + (3 - 2 a) x)/(3 - 2 c) ,
\\
f_{\sRD}(x,y) &= (2 + a - 5 c + (2 - a + c) x)/(4 - a - c) ,
\\
f_{\sUBR}(x,y) &=
(- 4 + 3 a + 2 b - c + (2 - a - b) x + (2 - b - c) y)/(a - c) ,
\\
f_{\sBR}(x,y) &=
( - 2 + 3 a + 2 b - 3 c + (2 - a - b) x + (3 - a - b) y)/(3 - 2 c) ,
\\
f_{\sBRD}(x,y) &= 
(a + 2 b - 5 c + (2 - a - b) x + (2 - b + c) y)/(4 - a - c).
\end{align*}
Then the graph of the function
 $z = f_{F}(x,y)$
is the perpendicular bisector of $\pU$ and $p_F$ 
for $F \in \BackRight \setminus \set{\sU}$.

Let 
\[ f_{\max}(x,y) = \max \set{f_F(x,y) \mid F \in \BackRight \setminus \set{\sU}} ,
\]
\[ R(F) = \set{(x,y) \in [a,1] \times [b,1] 
 \mid f_{F}(x,y) = f_{\max}(x,y) }
\]
for $F \in \BackRight \setminus \set{\sU}$.
Then a circumcenter $(x_0,y_0,z_0) = \rmc(F_1,F_2,F_3)$ is a Voronoi point
if and only if 
$(x_0,y_0) \in R(F_1) \cap R(F_2) \cap R(F_3)$
for $F_1,F_2,F_3 \in \BackRight \setminus\set{\sU}$. 
This is also equivalent to
\begin{equation}
\label{8}
 f_{\max}(x_0,y_0) = z_0 .
\end{equation}

Note that 
$R(\sUBR)$, $R(\sBR)$, $R(\sBRD)$ have empty interiors 
because
\begin{align*}
f_{\sUBR}(x,y) &< \max(f_{\sUB}(x,y), f_{\sUR}(x,y)),
\\
f_{\sBR}(x,y) &< \max(f_{\sB}(x,y), f_{\sUR}(x,y)),
\\
f_{\sBRD}(x,y) &< \max(f_{\sBD}(x,y), f_{\sUR}(x,y)),
\end{align*}
for $(x,y) \in [a,1)\times[b,1)$.

Denote by $v_0=(x_0,y_0,z_0)$.
If $p \in \Delta_{11}$, 
we can see that
$z_0 \ge f_{F}(x_0,y_0)$ for $v_0 = \rmc(\sB,\sD,\sR)$
with $F\in\set{\sUB, \sBD, \sRD, \sUR}$, 
so we have (\ref{8}) for $v_0 = \rmc(\sB,\sD,\sR)$.
We can also show that (\ref{8})
for 
$v_0 \in \set{\rmc(\sB,\sR,\sUR),  
\rmc(\sUB,\sB,\sUR)}$.
We also have (\ref{8}) 
for each $v_0 \in v(\sU)$ for other cases 
$p \in \Delta_{12}$, $\Delta_{21}$, $\Delta_{22A}$,
$\Delta_{22B}$, $\Delta_{31}$, $\Delta_{32A}$, 
$\Delta_{32B}$, $\Delta_{33}$.
\end{proof}

\begin{lemma}
Let $q \in \iota\Delta$ be a farthest point of $p \in \Delta$.
Then $q \in v(\sU)$.
\end{lemma}
\begin{proof}
If $p\in \Delta_{11}$,
we have
\[ \iota\Delta
 \subset t(\sfD,\sfR,\sfB) \cup t(\sR,\sUR,\sB) \cup 
t(\sUR,\sB,\sUB) .
\]
So $q$ is the circumcenter of either $t(\sfD,\sfR,\sfB)$,
$t(\sR,\sUR,\sB)$, $t(\sUR,\sB,\sUB)$.
Similar arguments hold for other cases.
\end{proof}

\section{The farthest point map}

The following theorem
shows the farthest point map on the $4$-cube.

\begin{theorem}
Let $q \in \iota\Delta$ be a farthest point of $p \in \Delta$.
If $p \in \Delta_{11}$, then
\begin{equation}
q = \rmc(\sB,\sD,\sR).
\label{7}
\end{equation}
If $p \in \Delta_{12}$, then
\[
 q = \rmc(\sBD,\sD,\sR).
\]
If $p \in \Delta_{21}$, then
\[
q \in  \set{ \rmc(\sB,\sD,\sRD), 
 \rmc(\sB,\sRD,\sR) } .
\]
If $p \in \Delta_{22A}$, then
\[
q \in  \set{ \rmc(\sBD,\sD,\sRD),
  \rmc(\sB,\sBD,\sRD), \rmc(\sB,\sRD,\sR) } .
\]
If $p \in \Delta_{22B}$, then
\[
q = \rmc(\sBD,\sD,\sRD) .
\]
If $p \in \Delta_{31}$, then
\[
q \in \set{ \rmc(\sB,\sD,\sRD), \rmc(\sB,\sRD,\sUR) } .
\]
If $p \in \Delta_{32A}$, then
\[
q \in \set{\rmc(\sBD,\sD,\sRD), 
 \rmc(\sB,\sBD,\sRD), \rmc(\sB,\sRD,\sUR) } .
\]
If $p \in \Delta_{32B}$, then
\[
q \in \set{\rmc(\sBD,\sD,\sRD), 
\rmc(\sBD,\sRD,\sUR), \rmc(\sB,\sBD,\sUR) } .
\]
If $p \in \Delta_{33}$, then
\[
q \in  \set{\rmc(\sBD,\sD,\sRD), 
\rmc(\sBD,\sRD,\sUR), 
\rmc(\sUB,\sBD,\sUR)} .
\]
\end{theorem}

\begin{proof}
If $p\in \Delta_{11}$, we have
\[ |\pU - \rmc(\sB,\sD,\sR)|
  > |\pU - \rmc(\sB,\sR,\sUR)|
  > |\pU - \rmc(\sUB,\sB,\sUR)| ,
\]
so the farthest point is $\rmc(\sB,\sD,\sR)$.
If $p\in \Delta_{12}$, we have
\[ |\pU - \rmc(\sBD,\sD,\sR)|
  > |\pU - \rmc(\sB,\sBD,\sR)|
  > |\pU - \rmc(\sB,\sR,\sUR)|
  > |\pU - \rmc(\sUB,\sB,\sUR)| ,
\]
so the farthest point is $\rmc(\sBD,\sD,\sR)$.
If $p\in \Delta_{21}$, 
\[ |\pU - \rmc(\sB, \sRD,\sR)|
  > |\pU - \rmc(\sB,\sR,\sUR)|
  > |\pU - \rmc(\sUB,\sB,\sUR)| , 
\]
so $\rmc(\sB,\sR,\sUR)$, $\rmc(\sUB,\sB,\sUR)$ are not farthest points.
If $p \in \Delta_{22A}$, then
\[
  |\pU - \rmc(\sB, \sRD,\sR)|
  > |\pU - \rmc(\sB, \sR,\sUR)|
  > |\pU - \rmc(\sUB,\sB,\sUR)| 
\]
so $\rmc(\sB, \sR,\sUR)$, $\rmc(\sUB, \sB, \sUR)$ are not farthest points.
If $p \in \Delta_{22B}$,
\begin{align*}
 & |\pU - \rmc(\sBD, \sD,\sRD)|
  > |\pU - \rmc(\sBD, \sRD,\sR)|
  > |\pU - \rmc(\sB,\sBD,\sR)|
\\
& > |\pU - \rmc(\sB, \sR,\sUR)|
  > |\pU - \rmc(\sUB,\sB,\sUR)| ,
\end{align*}
so the farthest point is $\rmc(\sBD,\sD,\sRD)$.  
If $p \in \Delta_{31} \cup \Delta_{32A}$, then
\[
  |\pU - \rmc(\sB, \sRD,\sUR)|
  > |\pU - \rmc(\sUB,\sB,\sUR)| ,
\]
so $\rmc(\sUB, \sB, \sUR)$ is not a farthest point.
If $p \in \Delta_{32B}$, then
\[
  |\pU - \rmc(\sB,\sBD,\sUR)|
  > |\pU - \rmc(\sUB,\sB,\sUR)| ,
\]
so $\rmc(\sUB, \sB, \sUR)$ is not a farthest point.
\end{proof}

It is not difficult to see that if $q$ is a farthest point of $p \in \Delta \cap H_{a=b}$, then $q \in \iota(\Delta \cap H_{a=b})$.
If $q$ is a farthest point of $p \in \Delta \cap H_{b=c}$, 
then $q \in  \iota(\Delta \cap H_{b=c})$.
If $q$ is a farthest point of $p \in \Delta \cap H_{c=0}$,
then $q \in \iota(\Delta \cap H_{c=0})$.
If $q \in \iota\Delta$ is a farthest point of
$p \in \Delta \cap H_{a=1/2}$,
then $r(q)$ is also a farthest point of $p$
where  $r = r_{H_{a=1/2}}$ is
the reflection with respect to $H_{a=1/2}$.

Next, we assume that 
$|p|$ is sufficiently small 
so that $p$ is close to a vertex of the $4$-cube. 
This will be used in the next section.
First, note that for a small neighborhood $U$ of the origin,
$U \cap \Delta_{12}$ and
$U \cap \Delta_{22B}$ are empty.

\begin{lemma}
\label{10}
Let $q \in \iota\Delta$ be a farthest point of $p\in\Delta$
and assume that $|p|$ is sufficiently small.
If $p \in \Delta_{11}$, then
\[ q = \rmc(\sB, \sD,\sR) .
\]
If $p \in \Delta_{21} \cup \Delta_{22A}$, 
\[ q = \rmc(\sB, \sRD, \sR) .
\]
If $p \in \Delta_{31} \cup \Delta_{32A}$, 
\[ q = \rmc(\sB, \sRD,\sUR) .
\]
If $p \in \Delta_{32B}$,
\[ q \in \set{ \rmc(\sBD, \sRD,\sUR), \rmc(\sB, \sBD, \sUR)} .
\]
If $p \in \Delta_{33}$, 
\[
q  \in \set{ \rmc(\sBD, \sRD,\sUR), \rmc(\sUB,\sBD,\sUR)} .
\]
\end{lemma}
\begin{proof}
Denote by $p = (a,b,c,0)$.
Let $p_{21}=(a_{21},b_{21},c_{21}):=(0.2864,0.09847,0.06212) \in \partial \Delta_{21}$
where we have $\psi_2(a_{21},c_{21})=0$, $\psi_1(b_{21},c_{21})=0$,  
$|\pU - \rmc(\sB,\sD,\sRD)| = |\pU-\rmc(\sB,\sRD,\sR)|$ at $p=p_{21}$.
For $p \in \Delta_{21}$ with $c<c_{21}$, we have
\[ |\pU - \rmc(\sB,\sD,\sRD)| < |\pU-\rmc(\sB,\sRD,\sR)|,
\]
so the farthest point is $\rmc(\sB,\sRD,\sR)$. 

Let $p_{22A1}=(a_{22A1},b_{22A1},c_{22A1})=(0.3282,0.3282,0.06898) \in \partial\Delta_{22A}$ 
where 
$a_{22A1}=b_{22A1}$,
$\psiur(a_{22A1},c_{22A1})=0$,
$|\pU - \rmc(\sBD,\sD,\sRD)| = |\pU - \rmc(\sB,\sBD,\sRD)|$
at $p=p_{22A1}$.
For $p \in \Delta_{22A}$ with $c<c_{22A1}$,
we have 
\[ |\pU - \rmc(\sBD,\sD,\sRD)| < |\pU - \rmc(\sB,\sBD,\sRD)| ,
\]
so $\rmc(\sBD,\sD,\sRD)$ is not a farthest point.

Let $p_{22A2} = (a_{22A2},b_{22A2},c_{22A2})=(0.2591,0.2591,0.05726)
 \in \Delta_{22A}$ where
$a_{22A2}=b_{22A2}$,
$\psiur(a_{22A2},c_{22A2})=0$,
$|\pU - \rmc(\sB,\sBD,\sRD)| = |\pU - \rmc(\sB,\sRD,\sR)|$
at $p=p_{22A2}$.
For $p \in \Delta_{22A}$ with $c<c_{22B2}$, 
we have
\[ |\pU - \rmc(\sB,\sBD,\sRD)| < |\pU - \rmc(\sB,\sRD,\sR)|,
\]
so $\rmc(\sB,\sBD,\sRD)$ is not a farthest point, and the farthest point is $\rmc(\sB,\sRD,\sR)$.


Let $p_{31}=(a_{31},b_{31},c_{31})=(0.2864, 0.09847,0.06212) \in \partial\Delta_{31}$ where
$\psiur(a_{31},c_{31})=0$,
$\psird(b_{31},c_{31})=0$, 
$|\pU - \rmc(\sB,\sD,\sRD)| = |\pU - \rmc(\sB,\sRD,\sUR)|$
at $p=p_{31}$.
For $p \in \Delta_{31}$ with $a<a_{31}$, we have 
\[ |\pU - \rmc(\sB,\sD,\sRD)| <|\pU - \rmc(\sB,\sRD,\sUR)|,
\]
so the farthest point is $\rmc(\sB,\sRD,\sUR)$.

Let $p_{32A1}=(a_{32A1},b_{32A1},c_{32A1})=(0.3282, 0.3282, 0.06898) \in \partial \Delta_{32A}$ where
$a_{32A1}=b_{32A1}$,
$\psiur(a_{32A1},c_{32A1})=0$, 
$|\pU-\rmc(\sBD,\sD,\sRD)| = |\pU-\rmc(\sB,\sBD,\sRD)|$
at $p=p_{32A1}$.
For $p\in\Delta_{32A}$ with $a<a_{32A1}$, we have
\[ |\pU-\rmc(\sBD,\sD,\sRD)| < |\pU-\rmc(\sB,\sBD,\sRD)|,
\]
so $\rmc(\sBD,\sD,\sRD)$ is not a farthest point.

Let $p_{32A2} = (a_{32A2},b_{32A2},c_{32A2})=(0.2591, 0.2591, 0.05726) \in \partial\Delta_{32A}$ where
$a_{32A2}=b_{32A2}$,
$\psiur(a_{32A2},c_{32A2})=0$, 
$|\pU-\rmc(\sB,\sBD,\sRD)| = |\pU-\rmc(\sB,\sRD,\sUR)|$
at $p=p_{32A2}$. 
For $p\in \Delta_{32A}$ with $a<a_{32A2}$,  
we have
\[ |\pU-\rmc(\sB,\sBD,\sRD)| < |\pU-\rmc(\sB,\sRD,\sUR)|,
\]
so the farthest point is $\rmc(\sB,\sRD,\sUR)$.

Let $p_{32B1}=(a_{32B1},b_{32B1},c_{32B1})=(0.3282, 0.3282,0.06898) \in \partial\Delta_{32B} \cap \partial\Delta_{33}$ where  
$a_{32B1}=b_{32B1}$,
$\psiur(a_{32B1},c_{32B1})=0$,
$\psiur(b_{32B1},c_{32B1})=0$,
$|\pU-\rmc(\sBD,\sD,\sRD)| = |\pU-\rmc(\sBD,\sRD,\sUR)|$ at
$p = p_{32B1}$.
For $p \in \Delta_{32B} \cup \Delta_{33}$ with $a<a_{32B1}$,
we have
\[ |\pU-\rmc(\sBD,\sD,\sRD)| < |\pU-\rmc(\sBD,\sRD,\sUR)|,
\]
so $\rmc(\sBD,\sD,\sRD)$ is not a farthest point of $p$.
\end{proof}

\section{Dynamics of the farthest point mapping}


In this section we show the following theorem, which was announced in
\cite{8OSME}.

\begin{theorem}
\label{17}
Let $q \in \iota\Delta$ be a farthest point of $p = (a,b,c,0) \in \Delta$.
Let $\set{q_j}_{j\ge1}$ be a sequence such that
$q_1 = \iota(q)$, and
$\iota(q_{j+1})$ is a farthest point of $q_j$, $j\ge1$.
Then there exists $0\le c' \le c$
such that
\[ \lim_{j\to\infty} q_j = (c',c',c',0).
\]
Furthermore, we have $c'>0$ if $c>0$.
\end{theorem}

First, we show that $\Delta_{11}$ is forward invariant under
the dynamics $\iota f$.

\begin{lemma}
Let $q \in \iota\Delta$ be a farthest point of $p =(a,b,c,0)\in \Delta_{11}$.
Then we have $\iota(q) \in \Delta_{11}$.
Define a sequence $\set{q_j}_{j\ge1}$
by $q_1 = \iota(q)$,
$q_{j+1}=\iota f(q_j)$, $j\ge1$.
Then
\[ \lim_{j\to\infty} q_j = (c,c,c,0).
\]
\end{lemma}
\begin{proof}
\[ \iota q = (\varphi_{1,c}(a), \varphi_{1,c}(b), c,0)
= \left( \frac{a +2c(1-c)}{3-2a}, \frac{b +2c(1-c)}{3-2b}, c, 0
 \right) .
\]
\end{proof}

If a forward orbit $\set{q_j}$ stays outside $\Delta_{11}$,
next Lemma~\ref{9} claims that
the limit must be the corner of the $4$-cube.

\begin{lemma}
\label{9}
Let $q\in\iota\Delta$ be a farthest point of 
$p \in \Delta \setminus \Delta_{11}$.
Let $\set{q_j}_{j\ge1}$ be a sequence such that
$q_1 = \iota(q)$,
$\iota(q_{j+1})$ is a farthest point of $q_j$, $j\ge1$.
Assume that $q_j \in \Delta \setminus \Delta_{11}$
for all $j\ge1$.
Then $\lim_{j\to\infty} q_j = 0$.
\end{lemma}
\begin{proof}
If $q' = (a',b',c',0) \in \Delta$ is a farthest point of 
$p=(a,b,c,0)\in\Delta \setminus \Delta_{11}$,
then $c < 0.102 < 1/6$, so
\begin{align*}
a' & \le \max\set{ \varphi_{1,c}(a), \varphi_{2,c}(a)}
 \le \frac{2}{3}(a-c)+c ,
\\
b' &\le \max\set{ \varphi_{1,c}(b), \varphi_{2,c}(b)} 
\le \frac{2}{3}(b-c)+c ,
\\
c' &\le c .
\end{align*}
\end{proof}

However, it can be shown that the assumption of Lemma~\ref{9}
 does not hold, which proves Theorem~\ref{17}.

\begin{lemma}
Let $q\in\iota\Delta$ be a farthest point of 
$p \in \Delta \setminus \Delta_{11}$.
Let $\set{q_j}_{j\ge1}$ be a sequence such that
$q_1 = \iota(q)$,
$\iota(q_{j+1})$ is a farthest point of $q_j$, $j\ge1$.
Then $q_j \in \Delta_{11}$ for sufficiently large $j$.
\end{lemma}
\begin{proof}
Suppose $p \in \Delta \setminus \Delta_{11}$.
By Lemma~\ref{9}, 
we can assume that $|p|$ is sufficiently small
so that Lemma~\ref{10} can be applied.

If $p \in \Delta_{21} \cup \Delta_{22}$, then
$q = \rmc(\sRD,\sR,\sB)$
and 
$\iota q \in \Delta_{11} \cup \Delta_{21} \cup \Delta_{22}$.
First we blow up the coordinates by considering the functions
\[ r_{xz}(a,b,c) = a/c,
\quad
r_{2}(a,b,c) = (a-c)/c^2.
\]
If $p=(a,b,c,0) \in \Delta_{21} \cup \Delta_{22}$ 
with $r_{xz}(p)>3$, we see that  
$r_{xz}(\iota q) < r_{xz}(p) - 1$.
So we have $r_{xz}(q_j) \le 3$ for some $j$.
If $p =(a,b,c,0) \in \Delta_{21} \cup \Delta_{22}$
with $r_{xz}(p) <3.535$,
we can show that $r_{xz}(\iota q) < r_{xz}(p)$, 
$r_{2}(p) > 6$,
and
\[ r_{2}(\iota q) < r_{2}(p) - 4.
\]
This implies that $q_j \in \Delta_{11}$ for some $j$.

If $p \in \Delta_{31} \cup \Delta_{32A}$,
then $q = \rmc(\sB,\sRD,\sUR)$ 
and $\iota q \in \Delta_{11} \cup \Delta_{21} \cup \Delta_{22} \cup \Delta_{31} \cup \Delta_{32A}$.
If $p \in \Delta_{31} \cup \Delta_{32A}$,
then $r_{xz}(p)= a/c >4$ and
\[ r_{xz}(\iota q) - r_{xz}(p) + 1 = - \frac{(a-3c)(a-c)}{2c(1-c)}<0.
\]
This implies that $q_j \in \Delta_{11} \cup \Delta_{21} \cup \Delta_{22}$
for some $j$.

If $p=(a,b,c,0)\in \Delta_{32B}$
and $a<0.3$,
then $q \in \set{ \rmc(\sBD,\sRD,\sUR), \rmc(\sB,\sBD,\sUR)}$
and $\iota q \in \Delta_{22} \cup \Delta_{32A}$.

If $p \in \Delta_{33}$, then 
$r_{yz}(p) > 4$ 
where $r_{yz}(a,b,c)=b/c$.
We have $q \in \set{\rmc(\sRD,\sUR,\sBD), \rmc(\sUR,\sBD,\sUB)}$.
so
\[ r_{yz}(q) - r_{yz}(p)+1 < 0,
\]
and $q_j \in \Delta \setminus \Delta_{33}$ for some $j$.
\end{proof}

\section{Intrinsic radius and diameter}

The intrinsic radius of a compact metric space $M$ is defined as
\[ \radius(M) := \inf \set{ d(p, f(p)) \mid p \in M },
\]
and the intrinsic diameter is defined as
\[ \diam(M) := \sup \set{ d(p, f(p)) \mid p \in M}
 = \sup \set{ d(p,q) \mid p,q \in M},
\]
where $f(p)$ denotes the farthest points of $p$.
It is easy to see that 
\[ 1/2 \le \radius(M)/\diam(M) \le 1.
\]

We have shown that
\begin{equation}
\label{18}
 \radius(\partial I^n) / \diam(\partial I^n) = 2 / \sqrt{n+2}
\end{equation}
for $n=2,3,4$.
On the $4$-cube, 
the intrinsic radius is given by
the facet centers $p=(1/2,1/2,1/2,0)$
and $q=(1/2,1/2,1/2,1) \in f(p)$, $d(p,f(p)) = 2$.
The intrinsic diameter is given by the corners $p=(0,0,0,0)$ and $q=(1,1,1,1)$,
where $d(p,q)=\sqrt{2^2 + 1^2 + 1^2} = \sqrt{6}$.

For $n\ge10$, the corner $q_0:=(1,\dots,1)$ 
is farther from $p_1:=(1/2,\dots,1/2,0)$ 
than the facet center $q_1:=(1/2,\dots,1/2,1)$,
\[ d(p_1,q_0)=\sqrt{(3/2)^2 + (n-2)\cdot 1^2} = \sqrt{n+7} / 2 
> d(p_1,q_1) = 2.
\]
It is conjectured for $n\ge10$ that 
the farthest point set $f(p)$ of any $p$ contains a corner,
the intrinsic radius $\radius(\partial I^n)$
is given by $p_1$ and $q_0$,
and the intrinsic diameter is given by the corners
$p_0=(0,\dots,0)$ and $q_0$, 
\[ d(p_0,q_0) = \sqrt{2^2+(n-2) \cdot 1^2}=\sqrt{n+2},
\]
so
\[ \radius(\partial I^n) / \diam(\partial I^n) = (\sqrt{n+7} / 2) /  \sqrt{n+2}
 \to 1/2
\] 
as $n\to \infty$.
It is expected that (\ref{18}) holds for $5 \le n \le 9$.

\section*{Acknowledgments}

The author would like to thank Sosuke Yagi, 
who worked on the farthest point map on the 4-cube
as part of his master's thesis.
This work was supported by JSPS Kakenhi Grant Number 24K06838.

\end{document}